\documentclass[12pt]{amsart}

\usepackage{amssymb,amsmath,amsthm,xcolor}
\usepackage{mathtools}
\usepackage{times}
\usepackage{microtype}
\usepackage[margin=1.2in]{geometry}
\usepackage{cite}
\usepackage[colorlinks=true,
 linkcolor=blue,citecolor=blue,urlcolor=blue]{hyperref}
\usepackage[nameinlink,capitalize,noabbrev]{cleveref}
\hypersetup{
 pdftitle={Polynomial maximal operators on Morrey sequence spaces},
 pdfauthor={Arash Ghorbanalizadeh and Zeinab Pourrezaei},
 pdfkeywords={Morrey sequence spaces; discrete polynomial maximal operators;
 polynomial averages; polynomial preimage estimates; pointwise convergence}
}

\numberwithin{equation}{section}
\allowdisplaybreaks

\newtheorem{theorem}{Theorem}[section]
\newtheorem{proposition}[theorem]{Proposition}
\newtheorem{lemma}[theorem]{Lemma}
\newtheorem{corollary}[theorem]{Corollary}
\theoremstyle{definition}
\newtheorem{definition}[theorem]{Definition}
\theoremstyle{remark}
\newtheorem{remark}[theorem]{Remark}

\newcommand{\Z}{\mathbb Z}
\newcommand{\N}{\mathbb N}
\newcommand{\NN}{\mathbb N_0}
\newcommand{\norm}[1]{\lVert #1\rVert}
\newcommand{\ind}{\mathbf 1}
\DeclareMathOperator{\sgn}{sgn}
\DeclareMathOperator{\supp}{supp}
\DeclareMathOperator{\dist}{dist}
\newcommand{\MP}{\mathcal M_P}
\newcommand{\AP}[1]{A_{#1}^{P}}
\newcommand{\Morrey}{m_p^\lambda(\Z)}
\newcommand{\mnorm}[1]{\lVert #1\rVert_{m_p^\lambda}}
\newcommand{\lpZ}{\ell^p(\Z)}
\newcommand{\linfZ}{\ell^\infty(\Z)}

\title[Polynomial maximal operators]
{Polynomial maximal operators on Morrey sequence spaces}
\date{}

\subjclass[2020]{Primary 42B25; Secondary 11B83, 46E30}\makeatletter
\def\@setsubjclass{%
  {\itshape\textup{2020} Mathematics Subject Classification.}\enspace
  \@subjclass\@addpunct.}
\makeatother
\keywords{Morrey sequence spaces, discrete polynomial maximal operators,
polynomial averages, polynomial preimage estimates, pointwise convergence}

\author[A. Ghorbanalizadeh]{Arash Ghorbanalizadeh}
\address{Department of Mathematics,\newline
Institute for Advanced Studies in Basic Sciences (IASBS),
Zanjan 45137-66731, Iran}
\email{ghorbanalizadeh@iasbs.ac.ir}
\email{gurbanalizade@gmail.com}
\thanks{Corresponding author: Arash Ghorbanalizadeh}

\author[Z. Pourrezaei]{Zeinab Pourrezaei}
\address{Department of Mathematics,\newline
Institute for Advanced Studies in Basic Sciences (IASBS),
Zanjan 45137-66731, Iran}
\email{z.pourrezaei@iasbs.ac.ir}

\begin{document}

\begin{abstract}
We study polynomial maximal operators on Morrey sequence spaces. For a nonconstant polynomial $P\in\mathbb Z[t]$ of degree $d$, we prove that the associated maximal operator is bounded on $m_p^\lambda(\mathbb Z)$ for $1< p < \infty$ and $0 \leq \lambda <1/d$. In the same range, the polynomial averages converge uniformly to zero at rate $O(N^{-(1-d\lambda)/p})$. For $d\geq2$, we show that the maximal operator is unbounded for $1/d\leq\lambda<1$ and construct sequences whose polynomial averages fail to converge at the origin for $1/d\leq\lambda\leq1$. At $\lambda=1$, the maximal operator is contractive for every $d\geq1$. For $d\geq2$, we also construct sequences for which the polynomial averages fail to converge pointwise when $1/d\leq\lambda\leq1$. The estimates use the difference between the number of averaging parameters and the length of the interval containing their polynomial  values. For the sharpness results, we use finite sets of polynomial values placed far apart. 
\end{abstract}
\maketitle

\section{Introduction and Preliminaries}\label{IP}

Polynomial averages are an important topic in discrete harmonic analysis
and ergodic theory.  Let $T$ be an invertible measure-preserving
transformation and let $P\in\Z[t]$.  A basic question in polynomial ergodic
theory is the convergence of averages of the form
\[
 \frac1N\sum_{k=1}^{N}F(T^{P(k)}x).
\]
Bourgain's work on pointwise ergodic theorems along polynomial sequences
showed the importance of maximal estimates in this problem
\cite{Bourgain1989}.  Later developments include discrete maximal and
variational inequalities \cite{MirekTrojan2016,MirekSteinTrojan2017,
MirekSteinTrojan2019,MirekSteinZorin2020}.  More recently, Krause,
Mirek, and Tao studied nonconventional bilinear polynomial ergodic
averages \cite{KrauseMirekTao2022}.  Their bilinear problem is different
from the linear maximal operator considered here, but it illustrates the
continuing role of polynomial averages in ergodic theory.

The translation model on $\Z$ leads to the polynomial maximal operator
studied in this paper.  Our starting point is the known estimate
\begin{equation}\label{eq:lp-input-intro}
\left\|
\sup_{N\geq1}\frac1N\sum_{m=1}^N
|g(\,\cdot-P(m))|
\right\|_{\lpZ}
\leq C_{p,P}\|g\|_{\lpZ},
\qquad 1<p<\infty.
\end{equation}
This follows from the polynomial-mapping maximal theorem in
\cite[Theorem~C]{MirekSteinTrojan2019}; see also
\cite{Bourgain1989,MirekTrojan2016}.  Polynomial averages on nilpotent
groups have also been studied in \cite{IonescuMagyarMirekSzarek2024}.

For $P(m)=m$, the averaging points fill an integer interval.  For a
polynomial of degree $d\geq2$, the values $P(m)$ with $|m|\leq N$
lie in an interval of length of order $N^d$, although they come from
only $2N+1$ parameters.  Moreover, each value has at most $d$ preimages.
Thus the number of sampled points and the size of the interval containing
them grow at different rates.  This difference is central to our results.

We ask what happens when the global $\ell^p$ condition in
\eqref{eq:lp-input-intro} is replaced by the local growth condition of a
Morrey sequence space.  Roughly speaking, this condition requires
\[
 \sum_{n\in J}|f(n)|^p\lesssim |J|^\lambda
\]
for every finite integer interval $J$.  Unlike $\ell^p$, Morrey spaces
with $0<\lambda<1$ can contain sequences whose total $\ell^p$ mass is
infinite.  For example, $h(n)=(1+|n|)^{-1/3}$ belongs to
$m_2^{1/3}(\Z)$ but not to $\ell^2(\Z)$.  Such slowly decaying sequences
show why it is useful to study the maximal operator beyond $\ell^p$.

Morrey sequence spaces and their inclusion properties were studied in
\cite{GunawanKikiantySchwanke2018}, and grand Morrey sequence spaces in
\cite{TulendinovaGhorbanalizadehZaighum2025}.  The discrete
Hardy--Littlewood maximal operator is bounded on Morrey sequence spaces
\cite[Theorem~3.2]{GunawanSchwanke2019}.  Our aim is to determine what
changes when ordinary integer shifts are replaced by polynomial shifts,
and to find the exact range of boundedness.

\subsection*{Notation and Morrey sequence spaces}
Throughout the paper, $\N=\{1,2,\ldots\}$ and $\NN=\{0,1,2,\ldots\}$.
For $a\in\Z$ and $r\in\NN$ we set
\begin{equation}\label{eq:interval}
I_a(r):=\{n\in\Z:\ |n-a|\leq r\},
\end{equation}
so that $|I_a(r)|=2r+1$.  The letter $C$ denotes a positive constant which may
change from line to line; its dependence on the parameters is indicated
whenever it is needed.  We write $A\lesssim B$ if $A\leq CB$, and we add the
parameters on which $C$ depends as subscripts, for example
$A\lesssim_{p,P}B$.

\begin{definition}\label[definition]{def:morrey}
Let $1\leq p<\infty$ and $0\leq\lambda\leq1$.  The Morrey sequence space
$m_p^\lambda=\Morrey$ consists of all complex sequences
$f=(f(n))_{n\in\Z}$ for which
\begin{equation}\label{eq:morrey-def}
\|f\|_{m_p^\lambda}
:=
\sup_{a\in\Z,\ r\in\NN}
\left(
(2r+1)^{-\lambda}
\sum_{n\in I_a(r)}|f(n)|^p
\right)^{1/p}
<\infty.
\end{equation}
\end{definition}

For $0\leq\lambda<1$ and $q=p/(1-\lambda)$, the norm
\eqref{eq:morrey-def} coincides with the norm of the discrete Morrey space
$\ell_q^p$ introduced in \cite{GunawanKikiantySchwanke2018}.  The two extreme
values of $\lambda$ give the familiar spaces
\begin{equation}\label{eq:end-spaces}
m_p^0(\Z)=\lpZ,
\qquad
m_p^1(\Z)=\linfZ,
\end{equation}
with equality of norms.  Indeed, the first identity follows by letting the
centered intervals $I_0(r)$ increase to $\Z$.  For the second one, every
local average of $|f|^p$ is bounded by $\|f\|_{\ell^\infty}^p$, while the
choice $r=0$ gives the reverse inequality.

We shall frequently use the Morrey condition on an arbitrary finite interval,
not necessarily centered.

\begin{lemma}\label[lemma]{lem:interval}
Let $1\leq p<\infty$, $0\leq\lambda\leq1$, and let $J\subset\Z$ be a
nonempty finite interval of integers. Then for every $f\in\Morrey$,
\begin{equation}\label{eq:arbitrary-interval}
\sum_{n\in J}|f(n)|^p\leq 2^\lambda|J|^\lambda\,\|f\|_{m_p^\lambda}^p .
\end{equation}
\end{lemma}

\begin{proof}
Write $J=\{u,u+1,\ldots,u+L-1\}$, where $L=|J|$, and put
$a:=u+\lfloor (L-1)/2\rfloor$ and $r:=\lceil (L-1)/2\rceil$.
Then $J\subset I_a(r)$ and $2r+1\leq L+1\leq 2L$.  Hence, by
\eqref{eq:morrey-def},
\[
\sum_{n\in J}|f(n)|^p
\leq \sum_{n\in I_a(r)}|f(n)|^p
\leq (2r+1)^\lambda\|f\|_{m_p^\lambda}^p
\leq (2L)^\lambda\|f\|_{m_p^\lambda}^p. \qedhere
\]
\end{proof}

\subsection*{Polynomial averages}
Let
\begin{equation}\label{eq:poly}
P(t)=a_dt^d+a_{d-1}t^{d-1}+\cdots+a_0\in\Z[t],
\qquad a_d\neq0,
\end{equation}
be a polynomial of degree $d\geq1$.  For $N\in\NN$ we define the polynomial
averages
\begin{equation}\label{eq:average}
\AP{N} f(n)
:=
\frac1{2N+1}
\sum_{|m|\leq N}|f(n-P(m))|,
\qquad n\in\Z,
\end{equation}
and the associated maximal operator
\begin{equation}\label{eq:maximal}
\MP f(n)
:=
\sup_{N\in\NN}\AP{N} f(n).
\end{equation}
Note that the absolute value is placed inside the average, so that $\MP$
dominates the maximal operator in which the absolute value is taken outside
the signed polynomial averages.

The Morrey norm \eqref{eq:morrey-def} is translation invariant, and
therefore the constant term $a_0$ plays no role in the boundedness
problem.  For $P(m)=m$, the operator is the discrete Hardy--Littlewood
maximal operator.  For $d\geq2$, the sparse polynomial shifts lead to
a different boundedness range.

\medskip
Our first main result proves that, for $1<p<\infty$, the operator
\[
\MP:\Morrey\longrightarrow\Morrey
\]
is bounded whenever $0\leq\lambda<1/d$.  The proof uses the known
polynomial maximal estimate on $\ell^p$ for the local part and a
local--far decomposition, as in classical Morrey-space arguments
(cf.~\cite{CF87,KRRS23,GunawanSchwanke2019}).  The far part is
controlled by counting the polynomial preimages of integer intervals.
This counting estimate is what gives the condition $\lambda<1/d$.

We then show that this condition is sharp when $d\geq2$.  For every
$1/d\leq\lambda<1$ we construct finitely supported sequences for which
$\MP$ is unbounded.  At $\lambda=1$, the space is $\ell^\infty$ and
$\MP$ is contractive.  Thus we obtain the complete boundedness range
for $1<p<\infty$.  We also prove that the polynomial averages converge
uniformly to zero at an explicit rate when $\lambda<1/d$.  In contrast,
for $d\geq2$ and $1/d\leq\lambda\leq1$, we construct sequences in
$\Morrey$ for which the averages fail to converge at the origin.  Hence
the same critical exponent $1/d$ governs both boundedness and the
convergence conclusion for all sequences in these spaces.

\medskip
The paper is organized as follows.  In \S~\ref{sec:MR} we state the main
results.  In \S~\ref{sec:DPME} we recall the discrete polynomial $\ell^p$
maximal estimate and prove the auxiliary lemmas, in particular the
polynomial preimage-counting lemma.  The main results are proved in
\S~\ref{sec:proofs}.

\section{Main results}\label{sec:MR}

The following results describe the boundedness of $\MP$ on Morrey sequence
spaces and the behavior of the polynomial averages.

\begin{theorem}\label{thm:main}
Let $P\in\Z[t]$ be nonconstant of degree $d\geq1$.  If $1<p<\infty$ and
$0\leq\lambda<1/d$, then there is a constant $C=C(p,\lambda,P)$ such that
\begin{equation}\label{eq:main}
 \norm{\MP f}_{m_p^\lambda}
 \leq C\mnorm{f}
\end{equation}
for every $f\in \Morrey$.
\end{theorem}

The next result gives a uniform bound for the polynomial averages.

\begin{proposition}\label[proposition]{prop:uniform-decay}
Let $P\in\Z[t]$ be nonconstant of degree $d\geq1$.  If $1<p<\infty$ and
$0\leq\lambda<1/d$, then there is a constant $C=C(p,\lambda,P)$ such that
\begin{equation}\label{eq:uniform-decay}
\sup_{n\in\Z} \AP{N} f(n)
\leq
C N^{-(1-d\lambda)/p}\,\mnorm{f},
\qquad N\geq1,
\end{equation}
for every $f\in \Morrey$.  In particular,
\[
\lim_{N\to\infty}\sup_{n\in\Z} \AP{N} f(n)=0.
\]
\end{proposition}

For instance, if $P(m)=m^d$, then \cref{thm:main,prop:uniform-decay} show that
the monomial maximal operator is bounded on $\Morrey$ and that the
corresponding averages converge uniformly to zero whenever $0\leq\lambda<1/d$.

For $d\geq2$, the range in \cref{thm:main} is sharp for $\lambda<1$.

\begin{theorem}\label{thm:sharpness}
Let $P\in\Z[t]$ be nonconstant of degree $d\geq2$ and let $1<p<\infty$.
For every $1/d\leq\lambda<1$, the operator
$\MP:\Morrey\longrightarrow\Morrey$ is unbounded.
More precisely, for every integer $K\geq1$ there is a finitely supported
sequence $f_K$ such that
\begin{equation}\label{eq:general-obstruction}
 \mnorm{f_K}\lesssim_{P,p}1,
 \qquad
 \mnorm{\MP f_K}\gtrsim_{p,\lambda}K^{(1-\lambda)/p}.
\end{equation}
\end{theorem}

\begin{proposition}\label[proposition]{prop:lambda-one}
For every nonconstant $P\in\Z[t]$ and every $1\leq p<\infty$,
\begin{equation}\label{eq:lambda-one}
 \norm{\MP f}_{m_p^1}\leq\norm{f}_{m_p^1}.
\end{equation}
\end{proposition}

Combining the three results above, we obtain the complete picture.

\begin{corollary}\label[corollary]{cor:exact-range}
Let $P\in\Z[t]$ have degree $d\geq2$ and let $1<p<\infty$.  Since $m_p^1(\Z)=\ell^\infty(\Z)$ with equality of norms,
the maximal operator $\MP$ is bounded
on $\Morrey$ for $0\leq\lambda<1/d$ and for $\lambda=1$, and it is unbounded
for every $1/d\leq\lambda<1$.
\end{corollary}

The last result concerns the averages themselves.  It is stated separately,
because the unboundedness of $\MP$ does not by itself imply the failure of
pointwise convergence.

\begin{proposition}\label[proposition]{prop:nonconvergence}
Let $P\in\Z[t]$ be nonconstant of degree $d\geq2$, let $1<p<\infty$, and let
$1/d\leq\lambda\leq1$.  Then there exists a nonnegative sequence
$f\in \Morrey$ such that
\begin{equation}\label{eq:nonconvergence}
 \limsup_{N\to\infty}\AP{N}f(0)>0,
 \qquad
 \liminf_{N\to\infty}\AP{N}f(0)=0.
\end{equation}
In particular, the polynomial averages need not converge pointwise on
$\Morrey$ in this range.
\end{proposition}

\section{Auxiliary estimates}\label{sec:DPME}

We begin with the discrete $\ell^p$ maximal estimate that will be used for the
local part of the proof of \cref{thm:main}.  The result is known in a much
more general form for polynomial mappings.  In the one-dimensional scalar
case, the required one-sided estimate is contained in
\cite[Theorem~C]{MirekSteinTrojan2019}; see also
\cite[Theorem~D]{MirekTrojan2016}.  We only need the following symmetric
consequence.

\begin{theorem}\label{thm:lp}
Let $P\in\Z[t]$ be nonconstant and let $1<p<\infty$.  Then there is a
constant $C_{p,P}>0$ such that
\begin{equation}\label{eq:lp-max}
 \norm{\MP g}_{\lpZ}\leq C_{p,P}\norm{g}_{\lpZ}
\end{equation}
for every $g\in\lpZ$.
\end{theorem}

\begin{proof}
First we explain how the cited theorem gives
\eqref{eq:lp-input-intro} when $P(0)$ is not zero.  Put
\[
 Q(t):=P(t)-P(0),\qquad h(n):=|g(n-P(0))|.
\]
Then $Q(0)=0$, $\|h\|_{\lpZ}=\|g\|_{\lpZ}$, and
$h(n-Q(m))=|g(n-P(m))|$.
Apply \cite[Theorem~C]{MirekSteinTrojan2019} to $Q$ with one nonzero
component, namely $h$.  This gives \eqref{eq:lp-input-intro}.

Now let $\widetilde P(t):=P(-t)$ and define
\[
M^{+}g:=\sup_{N\geq1}\frac1N\sum_{m=1}^{N}|g(\,\cdot-P(m))|,
\qquad
M^{-}g:=\sup_{N\geq1}\frac1N\sum_{m=1}^{N}|g(\,\cdot-\widetilde P(m))|.
\]
Applying \eqref{eq:lp-input-intro} to $P$ and $\widetilde P$, we obtain
$\|M^\pm g\|_{\lpZ}\leq C_{p,P}\|g\|_{\lpZ}$, after taking the larger of
the two constants.  For $N\geq1$, we separate the terms with $m=0$, $m>0$
and $m<0$.  This gives
\begin{align*}
\AP{N}g(n)
&\leq\frac{|g(n-P(0))|}{2N+1}+\frac{N}{2N+1}\bigl(M^{+}g(n)+M^{-}g(n)\bigr)\\
&\leq |g(n-P(0))|+\frac12\bigl(M^{+}g(n)+M^{-}g(n)\bigr),
\end{align*}
while $\AP{0}g(n)=|g(n-P(0))|$.  Hence
\[
\MP g\leq |g(\,\cdot-P(0))|+\frac12\bigl(M^{+}g+M^{-}g\bigr),
\]
and \eqref{eq:lp-max} follows by taking the $\lpZ$ norm.
\end{proof}

The next lemma controls the number of integers $m$ for which $P(m)$ lies in a
given integer interval.  It is elementary, but it is the main ingredient
of the Morrey estimate.

\begin{lemma}\label[lemma]{lem:counting}
Let $P$ be as in \eqref{eq:poly}.  Then there is a constant $C_d>0$,
depending only on $d$, such that for every finite integer interval $J$ with
$|J|=L\geq1$,
\begin{equation}\label{eq:counting}
 \#\{m\in\Z:\ P(m)\in J\}\leq C_d\,L^{1/d}.
\end{equation}
One may take $C_d=d\bigl(1+(2^{d-1}/d!)^{1/d}\bigr)$.
\end{lemma}

\begin{proof}
Let $S:=\{m\in\Z:\ P(m)\in J\}$ and $s:=\#S$.  Since $|P(m)|\to\infty$ as $|m|\to\infty$, the set $S$ is finite.
If $s\leq d$, then $s\leq dL^{1/d}$ and the result follows.
We therefore assume that $s\geq d+1$.

Let $x_1<x_2<\cdots<x_s$ be the elements of $S$, and put
\[
q:=\Bigl\lfloor\frac{s-1}{d}\Bigr\rfloor\ \ (\geq1),
\qquad
y_i:=x_{1+iq},\quad i=0,1,\ldots,d .
\]
The inequality $1+dq\leq s$ shows that all these points are defined.
The numbers $x_k$ are distinct integers in increasing order.  Hence
$x_v-x_u\geq v-u$ for $v>u$, and therefore
\begin{equation}\label{eq:separation}
 |y_i-y_j|\geq q\,|i-j|,\qquad 0\leq i,j\leq d.
\end{equation}
In particular, the points $y_i$ are distinct, and
\begin{equation}\label{eq:denom}
 \Bigl|\prod_{j\neq i}(y_i-y_j)\Bigr|
 \geq q^d\prod_{j\neq i}|i-j|=q^d\,i!\,(d-i)!\,.
\end{equation}

Next, write $J=\{u,u+1,\ldots,u+L-1\}$ and let $c:=u+\frac{L-1}{2}$ be its
midpoint.  Since $P(y_i)\in J$ for every $i$,
\begin{equation}\label{eq:values-small}
 |P(y_i)-c|\leq\frac{L-1}{2}\leq\frac L2,\qquad 0\leq i\leq d.
\end{equation}
Set $Q:=P-c$.  Since $Q$ has leading coefficient $a_d$, comparing
coefficients of $t^d$ in the Lagrange interpolation formula
\cite[\S~4.2, Eq.~(4.11)]{PhillipsTaylor1996} gives
\[
a_d=\sum_{i=0}^{d}
\frac{Q(y_i)}{\displaystyle\prod_{\substack{0\leq j\leq d\\j\ne i}}(y_i-y_j)}.
\]
By \eqref{eq:denom} and \eqref{eq:values-small},
\[
|a_d|
\leq\frac{L}{2q^d}\sum_{i=0}^{d}\frac1{i!\,(d-i)!}
=\frac{L}{2q^d}\cdot\frac{2^d}{d!}
=\frac{2^{d-1}}{d!}\cdot\frac{L}{q^d}\,,
\]
where we used $\sum_{i=0}^d\frac1{i!(d-i)!}=\frac1{d!}\sum_{i=0}^d\binom di
=\frac{2^d}{d!}$.  Since $P\in\Z[t]$ and $a_d\neq0$, we have $|a_d|\geq1$,
and therefore
\begin{equation}\label{eq:q-bound}
 q\leq\Bigl(\frac{2^{d-1}}{d!}\Bigr)^{1/d}L^{1/d}.
\end{equation}
Since $q=\lfloor(s-1)/d\rfloor$, we have $s-1<d(q+1)$.
Both $s$ and $d(q+1)$ are integers, so $s\leq d(q+1)$.
Together with \eqref{eq:q-bound} and $L^{1/d}\geq1$, this gives
\[
 s\leq d(q+1)\leq d\Bigl(\Bigl(\frac{2^{d-1}}{d!}\Bigr)^{1/d}+1\Bigr)L^{1/d}.
 \qedhere
\]
\end{proof}

\begin{remark}\label[remark]{rem:sharp-count}
The exponent $1/d$ in \eqref{eq:counting} cannot be improved.  For example, if
$P(m)=m^d$ and $J=\{0,1,\ldots,L-1\}$, then every integer
$0\leq m\leq (L-1)^{1/d}$ satisfies $P(m)\in J$, so the left-hand side of
\eqref{eq:counting} is of order $L^{1/d}$.
\end{remark}

We shall also use the following elementary growth estimate.

\begin{lemma}\label[lemma]{lem:growth}
Let $b_P:=\sum_{j=0}^d|a_j|$.  This is a fixed constant.  Then
\begin{equation}\label{eq:growth}
|P(m)|\leq b_P\,(1+|m|^d),\qquad m\in\Z.
\end{equation}
\end{lemma}

\begin{proof}
For $0\leq j\leq d$, we have $|m|^j\leq1+|m|^d$.
Thus
\[
 |P(m)|\leq\sum_{j=0}^d|a_j|\,|m|^j
 \leq b_P(1+|m|^d). \qedhere
\]
\end{proof}

The next lemma estimates one polynomial average on an interval $I$.
We will use it when $N^d$ is at least the number of points of $I$.

\begin{lemma}\label[lemma]{lem:one-scale}
Let $1\leq p<\infty$, $0\leq\lambda\leq1$, $a\in\Z$, $r\in\NN$, and put
$I:=I_a(r)$ and $L:=2r+1$.  If $N\in\N$ satisfies $N^d\geq L$, then
\begin{equation}\label{eq:one-scale}
 \sum_{n\in I}|\AP{N}g(n)|^p
 \leq C\,L^{1/d}N^{d\lambda-1}\,\mnorm{g}^p
\end{equation}
for every $g\in\Morrey$, where $C=C(d,\lambda,P)$.
\end{lemma}

\begin{proof}
By Jensen's inequality,
\[
|\AP{N}g(n)|^p\leq\frac1{2N+1}\sum_{|m|\leq N}|g(n-P(m))|^p .
\]
We sum over $n\in I$ and put $k=n-P(m)$.  For a fixed $k$, the
condition $n\in I$ becomes $P(m)\in I-k$.  Hence
\begin{equation}\label{eq:multiplicity-sum}
 \sum_{n\in I}|\AP{N}g(n)|^p
 \leq\frac1{2N+1}\sum_{k\in\Z}|g(k)|^p\,\nu(k),
\end{equation}
where
\[
 \nu(k):=\#\{m\in\Z:\ |m|\leq N,\ P(m)\in I-k\}.
\]
Since $I-k$ is an integer interval with exactly $L$ points,
\cref{lem:counting} gives
\begin{equation}\label{eq:multiplicity-bound}
 \nu(k)\leq C_d\,L^{1/d}\qquad\text{for all }k\in\Z .
\end{equation}
Moreover, $\nu(k)=0$ unless $k$ belongs to the interval
\[
 H:=\bigl[a-r-\Theta,\ a+r+\Theta\bigr]\cap\Z,
 \qquad
 \Theta:=\max_{|m|\leq N}|P(m)| .
\]
By \cref{lem:growth}, $\Theta\leq b_P(1+N^d)\leq2b_PN^d$.
Since $L\leq N^d$, we have
\[
|H|=L+2\Theta\leq(1+4b_P)N^d.
\]
Now we insert \eqref{eq:multiplicity-bound} in \eqref{eq:multiplicity-sum},
use $2N+1\geq N$, and apply \cref{lem:interval} on the interval $H$:
\[
 \sum_{n\in I}|\AP{N}g(n)|^p
 \leq\frac{C_dL^{1/d}}{N}\sum_{k\in H}|g(k)|^p
 \leq\frac{C_dL^{1/d}}{N}\,2^\lambda|H|^\lambda\mnorm{g}^p
 \leq C\,L^{1/d}N^{d\lambda-1}\mnorm{g}^p. \qedhere
\]
\end{proof}

\section{Proofs of the main results}\label{sec:proofs}

\subsection{Proof of \texorpdfstring{\cref{thm:main}}{Theorem 2.1}}

Fix an interval $I$.  We split $f$ into its restriction to a larger
interval $I^*$ and its restriction to the complement of $I^*$.
The number of points of $I^*$ is a fixed multiple of that of $I$.
We use \cref{thm:lp} for the first part and \cref{lem:one-scale} for the
second part.

\begin{proof}[Proof of \cref{thm:main}]
Fix $f\in\Morrey$ and an interval $I=I_a(r)$, and put $L:=2r+1$.
We will prove that
\begin{equation}\label{eq:local-goal}
 \sum_{n\in I}|\MP f(n)|^p\leq C\,L^\lambda\,\mnorm{f}^p,
\end{equation}
where $C=C(p,\lambda,P)$ does not depend on $a$, $r$ and $f$.

\medskip
\noindent\textit{Step 1: decomposition of $f$.}
Let $R:=\lceil L^{1/d}\rceil$.  Since $L\geq1$, we have
\begin{equation}\label{eq:R-comparable}
 L^{1/d}\leq R\leq2L^{1/d},\qquad R^d\leq2^dL .
\end{equation}
For $|m|\leq R$, \cref{lem:growth} gives
$|P(m)|\leq b_P(1+R^d)\leq b_P(1+2^d)L$.
Put $B_P:=\lceil b_P(1+2^d)\rceil$.  Then
\begin{equation}\label{eq:small-shift}
 |P(m)|\leq B_PL\qquad(|m|\leq R).
\end{equation}
Let $I^*:=I_a(r+B_PL)$.  We have $|I^*|=(1+2B_P)L$.
Also, $|n-P(m)-a|\leq r+B_PL$ when $n\in I$ and $|m|\leq R$.  Thus
\begin{equation}\label{eq:local-contained}
 n-P(m)\in I^*\qquad\text{whenever }n\in I,\ |m|\leq R.
\end{equation}
We write $f=f_1+f_2$, where $f_1:=f\ind_{I^*}$ and $f_2:=f\ind_{\Z\setminus I^*}$.
We have $\MP f\leq\MP f_1+\MP f_2$.
Since $(u+v)^p\leq2^{p-1}(u^p+v^p)$ for $u,v\geq0$, it is enough to
prove \eqref{eq:local-goal} for each part.

\medskip
\noindent\textit{Step 2: the local part.}
The sequence $f_1$ has finite support, so it belongs to $\lpZ$.  By
\cref{thm:lp} and \eqref{eq:morrey-def},
\begin{align}
 \sum_{n\in I}|\MP f_1(n)|^p
 &\leq\norm{\MP f_1}_{\lpZ}^p
  \leq C_{p,P}^p\,\norm{f_1}_{\lpZ}^p
  =C_{p,P}^p\sum_{k\in I^*}|f(k)|^p \notag\\
 &\leq C_{p,P}^p\,|I^*|^\lambda\mnorm{f}^p
  \leq C\,L^\lambda\mnorm{f}^p .
 \label{eq:local-bound}
\end{align}

\textit{Step 3: the averages of the far part.}
By \eqref{eq:local-contained} and the definition of $f_2$,
\begin{equation}\label{eq:far-small-zero}
 \AP{N}f_2(n)=0\qquad\text{for }n\in I,\ 0\leq N\leq R .
\end{equation}
It remains to consider $N>R$.  Put $N_j:=2^jR$ for $j\geq0$.
For each $N>R$, choose $j\geq1$ such that $N_{j-1}<N\leq N_j$.
Then $N_j<2N$.  All terms in the averages are nonnegative, so
\[
 \AP{N}f_2(n)\leq\frac{2N_j+1}{2N+1}\,\AP{N_j}f_2(n)<2\,\AP{N_j}f_2(n).
\]
Consequently, for $n\in I$,
\begin{equation}\label{eq:sup-dyadic}
 \MP f_2(n)\leq2\sup_{j\geq1}\AP{N_j}f_2(n),
 \qquad
 \bigl(\MP f_2(n)\bigr)^p\leq2^p\sum_{j\geq1}\bigl(\AP{N_j}f_2(n)\bigr)^p .
\end{equation}

\textit{Step 4: estimate of the far part.}
For every $j\geq1$ we have $N_j^d\geq R^d\geq L$, so \cref{lem:one-scale}
applies to $g=f_2$ and $N=N_j$.  Since $|f_2|\leq|f|$, we have
$\mnorm{f_2}\leq\mnorm{f}$, and hence
\[
 \sum_{n\in I}|\AP{N_j}f_2(n)|^p
 \leq C\,L^{1/d}N_j^{d\lambda-1}\mnorm{f}^p .
\]
Because $d\lambda-1<0$ and $R\geq L^{1/d}$,
\[
 L^{1/d}N_j^{d\lambda-1}
 =L^{1/d}\,2^{-j(1-d\lambda)}R^{d\lambda-1}
 \leq L^{1/d}\,2^{-j(1-d\lambda)}L^{(d\lambda-1)/d}
 =L^{\lambda}\,2^{-j(1-d\lambda)} .
\]
We sum \eqref{eq:sup-dyadic} over $n\in I$ and use the last two estimates.
This gives
\begin{equation}\label{eq:far-bound}
 \sum_{n\in I}|\MP f_2(n)|^p
 \leq C\,L^\lambda\mnorm{f}^p\sum_{j\geq1}2^{-j(1-d\lambda)}
 \leq C\,L^\lambda\mnorm{f}^p ,
\end{equation}
Here the series is finite because $1-d\lambda>0$; more precisely,
\[
 \sum_{j=1}^{\infty}2^{-j(1-d\lambda)}
 =\frac{1}{2^{1-d\lambda}-1}.
\]
This constant is allowed to depend on $\lambda$.

The estimates for $f_1$ and $f_2$ prove \eqref{eq:local-goal}.
Divide by $L^\lambda$, take the supremum over $a\in\Z$ and $r\in\NN$,
and then take the $p$th root.  We obtain \eqref{eq:main}.
\end{proof}

\subsection{Proof of \texorpdfstring{\cref{prop:uniform-decay}}{Proposition 2.2}}

\begin{proof}[Proof of \cref{prop:uniform-decay}]
Fix $f\in\Morrey$, $n\in\Z$ and $N\geq1$.  By Jensen's inequality,
\begin{equation}\label{eq:uniform-jensen}
 |\AP{N}f(n)|^p\leq\frac1{2N+1}\sum_{|m|\leq N}|f(n-P(m))|^p .
\end{equation}
Let $H_{n,N}$ be the smallest integer interval containing all points
$n-P(m)$ with $|m|\leq N$.  By \cref{lem:growth} and $N\geq1$,
\[
 |H_{n,N}|\leq1+2b_P(1+N^d)\leq(1+4b_P)N^d.
\]
For each $k$, the equation $P(m)=n-k$ has at most $d$ solutions.
Thus each term $|f(k)|^p$ occurs at most $d$ times in the sum in
\eqref{eq:uniform-jensen}.  By \cref{lem:interval},
\begin{equation}\label{eq:uniform-sum}
 \sum_{|m|\leq N}|f(n-P(m))|^p
 \leq d\sum_{k\in H_{n,N}}|f(k)|^p
 \leq d\,2^\lambda|H_{n,N}|^\lambda\mnorm{f}^p
 \leq C\,N^{d\lambda}\mnorm{f}^p .
\end{equation}
Using \eqref{eq:uniform-sum} in \eqref{eq:uniform-jensen}, we obtain
$|\AP{N}f(n)|^p\leq C N^{d\lambda-1}\mnorm{f}^p$.
Taking the $p$th root gives
\[
 \AP{N}f(n)\leq C\,N^{-(1-d\lambda)/p}\mnorm{f}.
\]
The constant is independent of $n$, so we may take the supremum over
$n\in\Z$.  This proves \eqref{eq:uniform-decay}.
Since $1-d\lambda>0$, the bound tends to zero as $N\to\infty$.
\end{proof}

\subsection{Sharpness and endpoint behavior}\label{sec:SEB}

We now prove unboundedness for $1/d\leq\lambda<1$.
For each integer $K\geq1$, we will construct a finitely supported sequence
whose Morrey norm stays bounded while the norm of its maximal function
tends to infinity as $K\to\infty$.

\begin{proof}[Proof of \cref{thm:sharpness}]
We construct a sparse input \(f_K\) whose maximal function is large at the consecutive integers \(1,\ldots,K\). To make the output large at \(j\), we set \(f_K=1\) at the points \(j-P(m)\) for \(R_j\leq m\leq 2R_j\). This makes at least a quarter of the terms in the average at \(j\) and scale \(2R_j\) equal to one. We choose the scales \(R_j\) so that the resulting sets are far enough apart that every interval of \(L\) integers contains at most \(C_P L^{1/d}\) points of the support. This keeps the input Morrey norm bounded independently of \(K\) when \(\lambda\geq 1/d\), while the \(K\) consecutive large output values force the output norm to grow when \(\lambda<1\).

Fix $1/d\leq\lambda<1$ and an integer $K\geq1$. 

\textit{Step 1: construction of $f_K$.}
Let $a_d$ be the leading coefficient of $P$ in \eqref{eq:poly}. Put $\sigma:=\sgn(a_d)$ and $Q:=\sigma P$.
Since $Q$ has positive leading coefficient, it is positive and strictly
increasing for all sufficiently large positive arguments.
Choose constants $0<c_Q\leq C_Q$ and an integer $m_0\geq1$, depending only
on $P$, such that $Q$ is strictly increasing for $m\geq m_0$ and
\begin{equation}\label{eq:polynomial-tail-growth}
 c_Q\,m^d\leq Q(m)\leq C_Q\,m^d,\qquad m\geq m_0 .
\end{equation}
Let $x_j:=j$ for $1\leq j\leq K$.  We choose integers $R_1<R_2<\cdots<R_K$ as
follows: $R_1\geq m_0$ is arbitrary, and, once $R_j$ has been chosen,
$R_{j+1}$ is taken so large that
\begin{equation}\label{eq:separated-scales}
 R_{j+1}\geq2R_j
 \qquad\text{and}\qquad
 c_Q\,R_{j+1}^d\geq4\bigl(C_Q(2R_j)^d+K\bigr).
\end{equation}
For $1\leq j\leq K$ we define
\begin{equation}\label{eq:general-blocks}
 E_j:=\{x_j-P(m):\ R_j\leq m\leq2R_j\},
 \qquad
 f_K:=\ind_{E_1\cup\cdots\cup E_K}.
\end{equation}
The values $P(m)$ are distinct for $m\geq m_0$, so $|E_j|=R_j+1$.
The scales $R_j$ may depend on $K$.  The constants $c_Q$, $C_Q$ and $m_0$
remain fixed.

\medskip
\noindent\textit{Step 2: the blocks are well separated.}
Let $1\leq i\leq j<K$, $u=x_i-P(m_i)\in E_i$ and
$v=x_{j+1}-P(m_{j+1})\in E_{j+1}$, where $R_i\leq m_i\leq2R_i$ and
$R_{j+1}\leq m_{j+1}\leq2R_{j+1}$.  Since $R_i\leq R_j$,
\eqref{eq:polynomial-tail-growth} gives
\[
 Q(m_i)\leq C_Q(2R_j)^d,\qquad Q(m_{j+1})\geq c_QR_{j+1}^d,
\]
and clearly $|x_i-x_{j+1}|\leq K$.  Because $P=\sigma Q$,
\begin{align*}
 |u-v|&=\bigl|x_i-x_{j+1}-\sigma\bigl(Q(m_i)-Q(m_{j+1})\bigr)\bigr|\\
 &\geq Q(m_{j+1})-Q(m_i)-K\\
 &\geq c_QR_{j+1}^d-C_Q(2R_j)^d-K
 \geq\tfrac34\,c_QR_{j+1}^d,
\end{align*}
where the last inequality follows from \eqref{eq:separated-scales}.  Hence
\begin{equation}\label{eq:block-separation}
 \dist\Bigl(E_{j+1},\ \bigcup_{i=1}^{j}E_i\Bigr)\geq\tfrac34\,c_QR_{j+1}^d .
\end{equation}
In particular, the sets $E_1,\ldots,E_K$ are pairwise disjoint.

\medskip
\noindent\textit{Step 3: counting the points of $\supp f_K$ in an interval.}
We claim that for every finite integer interval $J$ with $|J|=L$,
\begin{equation}\label{eq:general-count}
 \#\bigl(J\cap\supp f_K\bigr)\leq C_{P}\,L^{1/d} .
\end{equation}
If $J$ meets no block, the count is zero.
If it meets only $E_j$, apply \cref{lem:counting} to $x_j-J$ to obtain
$\#(J\cap E_j)\leq C_dL^{1/d}$.

Now suppose that $J$ meets at least two blocks, and let $j_0$ be the
largest index of a block meeting $J$.
The interval $J$ contains a point of $E_{j_0}$ and a point of an earlier
block.  Their distance is at most $L-1$.  By \eqref{eq:block-separation},
\[
 \frac34c_QR_{j_0}^d\leq L-1<L,
 \qquad
 R_{j_0}\leq\left(\frac{4}{3c_Q}\right)^{1/d}L^{1/d}.
\]
Also, $R_j\leq2^{-(j_0-j)}R_{j_0}$ for $j\leq j_0$.
It follows that
\[
 \#\bigl(J\cap\supp f_K\bigr)\leq\sum_{j=1}^{j_0}(R_j+1)
 \leq2\sum_{j=1}^{j_0}R_j\leq4R_{j_0}
 \leq4\left(\frac{4}{3c_Q}\right)^{1/d}L^{1/d}.
\]
Thus \eqref{eq:general-count} holds, for example, with
\[
 C_P:=\max\left\{C_d,\ 4\left(\frac{4}{3c_Q}\right)^{1/d}\right\}.
\]
In particular, this constant does not depend on $K$, $J$ or the chosen
scales.

\medskip
\noindent\textit{Step 4: the norm of $f_K$.}
Since $f_K$ is an indicator function, \eqref{eq:general-count} applied to
$J=I_a(r)$ gives
\[
 \sum_{n\in I_a(r)}|f_K(n)|^p=\#\bigl(I_a(r)\cap\supp f_K\bigr)
 \leq C_P(2r+1)^{1/d} .
\]
As $\lambda\geq1/d$, we obtain
\[
 \|f_K\|_{m_p^\lambda}^p
 \leq C_P\sup_{r\in\NN}(2r+1)^{1/d-\lambda}\leq C_P ,
\]
so that $\mnorm{f_K}\lesssim_{P,p}1$, with a bound independent of $K$.

\medskip
\noindent\textit{Step 5: the lower bound for $\MP f_K$.}
Fix $1\leq j\leq K$ and take $N=2R_j$.  For each of the $R_j+1$ integers
$m$ with $R_j\leq m\leq2R_j$ we have $x_j-P(m)\in E_j$, hence
$f_K(x_j-P(m))=1$.  Therefore
\begin{equation}\label{eq:general-output-pointwise}
 \MP f_K(x_j)\geq\AP{2R_j}f_K(x_j)\geq\frac{R_j+1}{4R_j+1}>\frac14 .
\end{equation}
The points $x_1,\ldots,x_K$ are consecutive integers.  Let
$r:=\lceil(K-1)/2\rceil$ and $I_K:=I_{r+1}(r)$.  Then
$\{1,\ldots,K\}\subset I_K$ and $|I_K|=2r+1\leq K+1$.  Using this interval in
the definition of the Morrey norm,
\[
 \|\MP f_K\|_{m_p^\lambda}^p
 \geq|I_K|^{-\lambda}\sum_{j=1}^{K}|\MP f_K(x_j)|^p
 \geq4^{-p}(K+1)^{-\lambda}K
 \geq4^{-p}2^{-\lambda}K^{1-\lambda}.
\]
Here we used $K+1\leq2K$.
Taking the $p$th root gives the second estimate in
\eqref{eq:general-obstruction}.

Since $\lambda<1$, the lower bound tends to infinity as $K\to\infty$.
The norms $\mnorm{f_K}$ remain bounded by a constant independent of $K$.
Therefore $\MP$ is unbounded on $\Morrey$.
\end{proof}

\begin{proof}[Proof of \cref{prop:lambda-one}]
Let $f\in m_p^1(\Z)=\linfZ$.  For every $n\in\Z$ and $N\in\NN$ we have
$\AP{N}f(n)\leq\|f\|_{\ell^\infty}$.  Taking the supremum over $N$ and $n$, and
using the equality of the norms of $m_p^1$ and $\ell^\infty$, we obtain
\eqref{eq:lambda-one}.
\end{proof}

\begin{proof}[Proof of \cref{cor:exact-range}]
The boundedness for $0\leq\lambda<1/d$ is \cref{thm:main}, the unboundedness
for $1/d\leq\lambda<1$ is \cref{thm:sharpness}, and the contractive estimate
at $\lambda=1$ is \cref{prop:lambda-one}.
\end{proof}

\subsection{Pointwise behavior of the polynomial averages}\label{sec:PBOPA}

We now construct a sequence whose averages do not converge at the origin.
Its support is the union of sets
$\{-P(k):A_j\leq k\leq2A_j\}$, where the integers $A_j$ increase rapidly.
The averages at $N=2A_j$ are bounded below by $1/4$.
The averages at $N=\lfloor A_j/4\rfloor$ tend to zero.

\begin{proof}[Proof of \cref{prop:nonconvergence}]
Fix $1/d\leq\lambda\leq1$.

\medskip
\noindent\textit{Step 1: choice of the parameters.}
Since $P(t)/(a_dt^d)\to1$ as $|t|\to\infty$, there is an integer $T\geq1$ such
that
\begin{equation}\label{eq:nonconv-growth}
 \tfrac12|a_d|\,|t|^d\leq|P(t)|\leq2|a_d|\,|t|^d\qquad(|t|\geq T).
\end{equation}
Let $B:=\max_{|m|<T}|P(m)|$.  We choose integers $A_1<A_2<\cdots$ such that
$A_1\geq\max\{T,4\}$, $B\leq\frac14|a_d|A_1^d$, and, for $j\geq2$,
\begin{equation}\label{eq:Aj}
 A_j\geq4A_{j-1}\qquad\text{and}\qquad
 A_j\geq j^2\sum_{i=1}^{j-1}(A_i+1).
\end{equation}
Let
\[
 S:=\bigcup_{j=1}^{\infty}\bigl([A_j,2A_j]\cap\N\bigr),
 \qquad
 f:=\ind_{\{-P(k):\ k\in S\}} .
\]
By the first condition in \eqref{eq:Aj}, the blocks $[A_j,2A_j]$ are pairwise
disjoint, and $f$ takes only the values $0$ and $1$.

\medskip
\noindent\textit{Step 2: $f\in\Morrey$.}
Let $J$ be a finite integer interval with $|J|=L$.  Since $-J$ is also an
interval with $L$ points, \cref{lem:counting} gives
\[
 \#(J\cap\supp f)\leq\#\{k\in\Z:\ P(k)\in-J\}\leq C_dL^{1/d} .
\]
Hence, for all $a\in\Z$ and $r\in\NN$, $\sum_{n\in I_a(r)}|f(n)|^p\leq
C_d(2r+1)^{1/d}$, and, since $\lambda\geq1/d$,
\[
 (2r+1)^{-\lambda}\sum_{n\in I_a(r)}|f(n)|^p\leq C_d(2r+1)^{1/d-\lambda}\leq C_d .
\]
Thus $f\in\Morrey$.

\medskip
\noindent\textit{Step 3: the limsup is positive.}
Take $N=2A_j$.  Every integer $m$ with $A_j\leq m\leq2A_j$ belongs to $S$, so
$f(-P(m))=1$ for these $A_j+1$ values of $m$.  Hence
\[
 \AP{2A_j}f(0)=\frac1{4A_j+1}\sum_{|m|\leq2A_j}f(-P(m))
 \geq\frac{A_j+1}{4A_j+1}>\frac14 ,
\]
and so $\limsup_{N\to\infty}\AP{N}f(0)\geq\frac14$.

\medskip
\noindent\textit{Step 4: the liminf is zero.}
Let $N_j:=\lfloor A_j/4\rfloor$.  Since $A_j/4\geq1$, we have
$N_j\geq A_j/8$; in particular $N_j\to\infty$.  Let $|m|\leq N_j$ be such that
$f(-P(m))=1$.  By the definition of $f$ there is $k\in S$ with
\begin{equation}\label{eq:collision}
 P(m)=P(k).
\end{equation}
We show that $k<A_j$.
For every $k\geq A_j$, \eqref{eq:nonconv-growth} gives
$|P(k)|\geq\frac12|a_d|A_j^d$.
For $|m|\leq N_j\leq A_j/4$, there are two cases.
If $|m|<T$, then
$|P(m)|\leq B\leq\frac14|a_d|A_j^d$.
If $|m|\geq T$, then, since $d\geq2$,
\[
 |P(m)|\leq2|a_d|\left(\frac{A_j}{4}\right)^d
 \leq\frac18|a_d|A_j^d.
\]
Thus $P(m)$ cannot equal $P(k)$ for any $k\geq A_j$.
The equality \eqref{eq:collision} therefore implies $k<A_j$.
This argument also covers negative values of $m$.

Consequently, $k$ belongs to one of the blocks $[A_i,2A_i]$ with $i<j$.
There are $\sum_{i<j}(A_i+1)$ possible values of $k$.
For each of them, $P(m)=P(k)$ has at most $d$ solutions.
For $j\geq2$, we can now use \eqref{eq:Aj} and
$2N_j+1\geq N_j\geq A_j/8$ to obtain
\[
 \AP{N_j}f(0)=\frac1{2N_j+1}\sum_{|m|\leq N_j}f(-P(m))
 \leq\frac{d\sum_{i=1}^{j-1}(A_i+1)}{2N_j+1}
 \leq\frac{8d}{A_j}\sum_{i=1}^{j-1}(A_i+1)\leq\frac{8d}{j^2}\longrightarrow0 .
\]
As the averages are nonnegative and $N_j\to\infty$, we conclude that
$\liminf_{N\to\infty}\AP{N}f(0)=0$.

We have proved that the limsup is at least $1/4$ and the liminf is zero.
This gives \eqref{eq:nonconvergence}.
\end{proof}

\begin{remark}
For $d\geq2$ and $1<p<\infty$, the two propositions give the following
complete statement.  If $0\leq\lambda<1/d$, then
$\sup_{n\in\Z}\AP{N}f(n)\to0$ for every $f\in\Morrey$.
If $1/d\leq\lambda\leq1$, there is a sequence in $\Morrey$ whose polynomial
averages do not converge at the origin.
\end{remark}

\section*{Statements and Declarations}

\noindent\textbf{Funding.}
The authors received no funding for this work.

\noindent\textbf{Competing interests.}
The authors declare that they have no competing interests.

\noindent\textbf{Use of generative AI.}
During the preparation of this manuscript, the authors used ChatGPT for
language editing and to improve the clarity of the presentation.


\begin{thebibliography}{99}

\bibitem{Bourgain1989}
J.~Bourgain,
\newblock Pointwise ergodic theorems for arithmetic sets,
\newblock \emph{Publ. Math. Inst. Hautes \'{E}tudes Sci.} \textbf{69} (1989),
5--41.

\bibitem{CF87}
F.~Chiarenza and M.~Frasca,
\newblock Morrey spaces and Hardy--Littlewood maximal function,
\newblock \emph{Rend. Mat. Appl. (7)} \textbf{7} (1987), no.~3--4, 273--279.

\bibitem{GunawanKikiantySchwanke2018}
H.~Gunawan, E.~Kikianty, and C.~Schwanke,
\newblock Discrete Morrey spaces and their inclusion properties,
\newblock \emph{Math. Nachr.} \textbf{291} (2018), 1283--1296.


\bibitem{GunawanSchwanke2019}
H.~Gunawan and C.~Schwanke,
\newblock The Hardy--Littlewood maximal operator on discrete Morrey spaces,
\newblock \emph{Mediterr. J. Math.} \textbf{16} (2019), Paper No.~24, 12 pp.


\bibitem{IonescuMagyarMirekSzarek2024}
A.~D. Ionescu, A.~Magyar, M.~Mirek, and T.~Z. Szarek,
\newblock Polynomial averages and pointwise ergodic theorems on nilpotent
groups,
\newblock \emph{Invent. Math.} \textbf{231} (2023), 1023--1140.


\bibitem{KRRS23}
A.~Kassymov, M.~A. Ragusa, M.~Ruzhansky, and D.~Suragan,
\newblock Stein--Weiss--Adams inequality on Morrey spaces,
\newblock \emph{J. Funct. Anal.} \textbf{285} (2023), no.~11,
Paper No.~110152, 30 pp.


\bibitem{KrauseMirekTao2022}
B.~Krause, M.~Mirek, and T.~Tao,
\newblock Pointwise ergodic theorems for non-conventional bilinear polynomial
averages,
\newblock \emph{Ann. of Math. (2)} \textbf{195} (2022), 997--1109.


\bibitem{LukyanovTerras2003}
S.~Lukyanov and V.~Terras,
\newblock Long-distance asymptotics of spin--spin correlation functions for the XXZ spin chain,
\newblock \emph{Nucl. Phys. B} \textbf{654} (2003), no.~3, 323--356.


\bibitem{MirekSteinTrojan2017}
M.~Mirek, E.~M. Stein, and B.~Trojan,
\newblock $\ell^p(\Z^d)$-estimates for discrete operators of Radon type:
variational estimates,
\newblock \emph{Invent. Math.} \textbf{209} (2017), 665--748.


\bibitem{MirekSteinTrojan2019}
M.~Mirek, E.~M. Stein, and B.~Trojan,
\newblock $\ell^p(\Z^d)$-estimates for discrete operators of Radon type:
maximal functions and vector-valued estimates,
\newblock \emph{J. Funct. Anal.} \textbf{277} (2019), 2471--2521.


\bibitem{MirekSteinZorin2020}
M.~Mirek, E.~M. Stein, and P.~Zorin-Kranich,
\newblock Jump inequalities for translation-invariant operators of Radon type
on $\Z^d$,
\newblock \emph{Adv. Math.} \textbf{365} (2020), 107065.


\bibitem{MirekTrojan2016}
M.~Mirek and B.~Trojan,
\newblock Discrete maximal functions in higher dimensions and applications to
ergodic theory,
\newblock \emph{Amer. J. Math.} \textbf{138} (2016), 1495--1532.


\bibitem{TulendinovaGhorbanalizadehZaighum2025}
A.~Tulendinova, A.~Ghorbanalizadeh, and M.~A. Zaighum,
\newblock On grand Morrey sequence spaces,
\newblock \emph{J. Pseudo-Differ. Oper. Appl.} \textbf{16} (2025),
Paper No.~73.


\bibitem{PhillipsTaylor1996}
G.~M. Phillips and P.~J. Taylor,
\newblock \emph{Theory and Applications of Numerical Analysis},
\newblock 2nd ed., Academic Press, 1996.

\end{thebibliography}
\end{document}